\documentstyle[titlepage]{article}
\title{Generalized Futaki Invariant\\ of Almost Fano Toric Varieties,\\Examples}
\author{M.Yotov\thanks{Partially supported by University of Oslo, and by a grant from EPDI.}\\$ $\\
author's permanent address:\\
University of Sofia\\
Department of Mathematics and Informatics\\
Chair of Geometry\\
5 ``Games Bourchier''Blvd.\\
1126 Sofia, Bulgaria.}
\begin{document}
\maketitle
\begin{abstract}
The interpretation, due to T. Mabuchi, of the classical Futaki invariant of Fano toric manifolds is extended to the case of the Generalized Futaki invariant, introduced by W.Ding and G.Tian, of almost Fano toric varieties. As an application it is shown that the real part of the Generalized Futaki invariant is positive for all degenerations of the Fano manifold $V_{38}$\,, obtained by intersection of the Veronese embedding of ${\bf P}^3\times{\bf P}^2 \subset {\bf P}^{11}$\,\, with codimension-two hyperplanes.
\end{abstract}

\newtheorem{theo}{Theorem}[section]
\newtheorem{prop}{Proposition}[theo]
\newtheorem{lem}{Lemma}[theo]
\newtheorem{cor}{Corollary}[theo]
\newtheorem{rem}{Remark}[section]
\newenvironment{definbas}{\begin{defin}\em}{\end{defin}}
\newtheorem{defin}{Definition}[section]

\section{Introduction}

The aim of this paper is to present the result of Toshiki Mabuchi on the calculation of the Futaki invariant of toric Fano manifolds, but adapted to the case of almost Fano toric varieties. The reason of doing this will be explained in what follows.

The famous Futaki invariant, when non-zero, serves as an obstruction to the existence of Einstein-K\"ahler metrics on a Fano manifold. Being a character of the Lie algebra of holomorphic vector fields, it is identically zero and gives no information in the case when the manifold has a semi-simple group of holomorphic automorphisms. Apart from this case, the calculation of the Futaki invariant in general is not easy. There is one beautiful case of Fano manifolds where this invariant can be calculated completely and given a nice geometric interpretaton - this is the case of toric Fano varieties ( which is due to Toshiki Mabuchi \cite{M1}). It is well known (see \cite{Fut1}, Chapter III) that in general the Futaki invariant can be interpreted as the baricentre of certain moment map of the manifold. But in the case of toric Fano manifolds Mabuchi has proved that it can be viewed as the baricentre, with respect to certain Duistermaat-Heckman measure, of the polytope corresponding to the anticanonical line bundle of the manifold.

To overcome the inconvenience of the case of semi-simple (in particular - finite) groups, Ding and Tian have introduced the Generalized Futaki invariant for almost Fano varieties and replaced the existence problem on the initial Fano manifold with the properties of some its degenerations. To be more precise, they consider degenerations of the Fano manifold which have almost Fano varieties as singular fibres and have naturally defined vector fields which can be restricted to those fibres. Then, nonpositivity of the real part of the Generalized Futaki invariant of some degenerated fibre, evaluated on the corresponding restricted vector field, serves as an obstruction to the existence of Einstein-K\"ahler metrics on the initial Fano manifold.

The reason of extending Mabuchi's to the case of almost Fano toric varieties is given by the fact that there are several examples of Fano manifolds which can be degenerated to almost Fano toric varieties so that the theory of Ding and Tian applies.

The structure of the paper is as follows. In Section 2 we present the definition of the generalized Futaki invariant according to Ding and Tian (\cite{D-T}) with a slight modification. In particular, we show that it depends on the integration over open and dense subsets of the variety only. This is the crucial observation for the rest of the note - it allows us to interpret the Generalized Futaki invariant of almost Fano toric varieties in the same way as Mabuchi does. We remind the definition of Weakly K-stability (due to G.Tian \cite{Tian}), and formulate Tian's theorem which shows the connection between Weakly K-stability and existence of Einstein-K\"ahler metrics on Fano manifolds. Section 3 presents the basics of the Moment-Map theory and introduces the moment map of ${\bf P}^N$\,used in the paper. Section 4 presents the basics of the theory of toric varieties in the scope we need it, and gives the description of the Generalized Futaki Invariant via the corresponding moment map. In particular, Theorem \ref{torus} gives an explicit formula how to calculate the real part of the Generalized Futaki Invariant of almost Fano toric varieties. In Section 5 we study the special degenerations of the Fano manifold $V_{38}$\,- the blowup of ${\bf P}^3$\, in a twisted cubic curve, obtained by intersection of the Veronese embedding of ${\bf P}^3\times{\bf P}^2 \,\subset {\bf P}^{11}$\,\, by codimension-two hyperplanes. We show that there are only three types of such degenerations with toric singular fibres. For all of them the real part of the Generalized Futaki invariant is positive (when evaluated on the vector fields defined by the degenerations).

This paper was prepared during the author's stay as a Post-Doc Fellow of the European Post-Doctoral Institute for Mathematical Sciences at the University of Oslo. He is deeply grateful to EPDI for giving him the opportunity to work in the very stimulating atmosphere of the Algebraic Geometry seminars of this University! Special thanks are due to K.Ranestad for the fruitful discussions the author has had with him, and for the permanent encouragement during the preparation of this paper.

\section{Almost Fano varieties and the Generalized Futaki invariant}

Let us begin with the definition of the classical Futaki invariant.

Suppose we are given a Fano manifold $X$. By definition this means that $X$ is compact complex manifold with positive first Chern class. By the famous Kodaira theorem we can choose a K\"ahler metric $g$ on $X$, which K\"ahler form $\omega _g$ represents the first Chern class of $X$ : $\omega _g \in 2\pi c_1(X)$. In local coordinates $z_1,\ldots ,z_n$ the K\"ahler form of $g$ can be expressed as

$$\omega _g = \sqrt {-1}.g_{i\bar j}dz^i\wedge d\bar{z^j}.$$

By definition, the Ricci form of $g$ is a (1,1)-form locally expressed as

$$ Ric(g) := - \sqrt{-1}.\partial \bar{\partial}\log (\det g_{i\bar{j}}).$$

The metric $g$ defines a volume form $dV_g := (1/n!)\omega _g^{\wedge n}.$

It is well known that the Ricci curvature of any K\"ahler metric represents the first Chern class of the manifold, up to a multiplication by $2\pi$. So, there exists a smooth function $f_g$ on $X$ such that

$$Ric(g) - \omega _g = \sqrt{-1}\partial \bar{\partial}f_g.$$

Now we define a complex linear map from the Lie algebra of holomorphic vector fields of $X$ to the complex numbers

$$F : Lie(Aut(X)) \longrightarrow {\bf C}$$
by assigning to a vector field $\xi$ the complex number $F(\xi) := \int_{X}\xi(f_g)dV_g$.

\begin{theo}[Futaki]
1) The map $F$ does not depend on the choice of the metric $g$ representing the first Chern class of $X$;

2) The map $F$ is a Lie-algebra homomorphism;

3) If $X$ is Einstein-K\"ahler, i.e., for some $g$ the corresponding $f_g$ is a constant function, then $F$ is identically zero.

\end{theo}

Ding and Tian \cite{D-T} have generalized the notion of Futaki invariant to the case of normal almost Fano varieties. The general setup goes as follows.

A compact complex variety $X$ is called almost Fano, if it is normal and there exists an ample line bundle on $L$ it, such that the restriction of ${\cal L}^*$ to the regular part of $X$ is isomorphic to some positive tensor power of the dualizing sheaf of $X$ :

$$ {\cal L}^* \cong \omega _{X_{reg}}^k.$$

In order to define the Generalized Futaki invariant we have to generalize the notions of a metric and holomorphic vector fields on an almost Fano variety.

\begin{itemize}
\item Admissible K\"ahler metrics on $X$ (with respect to $\omega _X^*$).
\end{itemize}

Let $L^m$ be very ample and $\phi _m : X \rightarrow {\bf P}^N$ be the corresponding embedding. If $\omega$ is a K\"ahler metric on ${\bf P}^N$, then $\phi _m^*\omega$ (on $X_{reg}$) is called a K\"ahler metric on $X$. The most important K\"ahler metrics on $X$ for us are of the form $\omega = {1\over k.m}\phi _m^*\omega '$, where $\omega '$ represents the first Chern class of the hyperplane line bundle $[H]$ over $X$, multiplied by $2\pi $.

{\bf Fact.} There is a smooth metric $\tilde{h}$ on $[H]$ such that $$\omega ' = Ric (\tilde{h}) := - \sqrt {-1}.\partial \bar \partial \log (\tilde{h})$$.

Hence, all our considerations will depend on a metric $\tilde{h}$ on $[H]$ for which $Ric(\tilde{h}) > 0$.

\begin{itemize}
\item Admissible vector fields on $X$ (with respect to $\omega _X^*$).
\end{itemize}
Let $\phi _m : X \rightarrow {\bf P}^N$ be as above, and let $\tilde{\xi} \in \Gamma ({\bf P}^N, {\cal T_{{\bf P}^N}})$ be a holomorphic vector field on the ambient manifold which is tangent to $\phi _m(X)$ along $\phi _m(X_{reg})$. Then the restriction $\xi := \tilde{\xi }_{\mid \phi _m(X_{reg})}$ is called an admissible vector field on $X$.

{\bf Facts.}
\begin{enumerate}
\item For an almost Fano variety $X$ the group of holomorphic automorphisms  $G := Aut(X)$ is a linear algebraic group. More precisely, $G = Aut(X, {\bf P}^N)$, where the righthand side denotes the group of holomorphic automorphisms of ${\bf P}^N$ leaving $\phi _m(X)$ fixed.

\item Let $\eta (X)$ be the set of admissible vector fields on $X$. Then, $\eta (X) = Lie G$, where the righthand side denotes the Lie algebra of the group $G$. Hence, the admissible vector fields are determined intrinsically and form a Lie algebra.
\end{enumerate}

The first assertion above follows from the fact that, since $X$\, is normal, $G$ acts naturally on $L^{\otimes m}$ lifting its action on $X$. The second one follows from the fact that the holomorphic vector fields on ${\bf P}^N$ which are tangent to $X$ along $X_{reg}$ are in bijection with the one-parametric subgroups of $Aut({\bf P}^N)$ which lie in $G \subset Aut({\bf P}^N)$.
\par
Let $\pi : \tilde{X} \rightarrow X$ be any ($G$ - equivariant) desingularization of $X$, and let $E$ be the exceptional subscheme of $\pi $. Let further, $K_{\tilde{X}}$ be the canonical line bundle of $\tilde{X}$. The line bundle $K_{\tilde{X}}^{\otimes k} \otimes \pi ^*L$ is trivial on $\tilde{X} - E$ and has a section over $\tilde{X}$ defined by $\pi$. Hence, it is a line bundle $[\tilde{D}]$, corresponding to an  effective divisor $\tilde{D}$ with support in $Supp(E)$.

By taking the k-th square root of the metric on $h := \phi _m^*(\tilde{h})$ on $L^{\otimes m}$ we get a metric $l$ on $L$, such that on $X_{reg}$ we have $Ric(l) = k.\omega$.

If we set $\tilde{\omega} := (\phi _m \circ \pi)^*\omega '$, then for the metric $\pi ^*(l)$ on $\pi ^*(L)$ we have  $Ric(\pi ^*(l)) = k.\tilde{\omega}$.
\par

Denote by $t$ any smooth metric on $K_{\tilde{X}}^*$. Then we have an induced smooth metric $\tilde{d} := \pi ^*l \otimes t^{\otimes(-k)}$ on the line bundle $[\tilde{D}]$.

The K\"ahler form $\tilde{\omega}$ defines a volume form $dV_{\tilde{\omega}} := (1/n!).\tilde{\omega}^{\wedge n}$ which is smooth on $\tilde{X}$\,, and is zero only on $Supp E$. On the other hand $t$, being a metric on $K_{\tilde{X}}^*$, can be regarded as a positive volume form $dV_t$ on $\tilde{X}$, too. Hence there is a well defined function $\tilde{\omega}/t$ on $\tilde{X}$, such that $dV_{\tilde{X}} = (\tilde{\omega}/t).dV_t$.

In these notations we have that on $\tilde{X}\,-\,E$

$$ Ric(\tilde{\omega}) - \tilde{\omega} = (-1/k).Ric(\tilde{d}.(t/\tilde{\omega})^k).$$

Since $[D]$ is trivial on $\tilde{X}\,-\,E$, then there exists a smooth positive function $d$ on $\tilde{X}\,-\,E$ which represents $\tilde{d}$ there. Define on $\tilde{X}\,-\,E$ the function

$$\varphi := (1/k).\log (d.(t/\tilde{\omega})^k).$$  
This function is uniquely determined, up to an additive scalar, by the form $\tilde{\omega}$, hence by the metric $\tilde{h}$, and we have on $\tilde{X}\,-\,E$

$$ Ric(\tilde{\omega}) - \tilde{\omega} = \sqrt{-1}\partial\bar{\partial} \varphi.$$
In particular, $\varphi $ does not depend on the desingularization $\pi$ of $X$.

\begin{itemize}
\item
{\bf Important fact.}
If $\tilde{\xi}$ is a smooth vector field on $\tilde{X}\,-\,E$ which is continuous on the whole $\tilde{X}$, then

$$\int_{\tilde{X}\,-\,E}\,\tilde{\xi}(\varphi)\,dV_{\tilde{\omega}}$$
is a well defined finite number. In particular, for any admissible vector field $\xi \in \eta (X)$ we can define

$$F_{\varphi}(\xi) := \int_{\tilde{X}\,-\,E}\tilde{\xi}(\varphi)\,dV_{\tilde{\omega}},$$
where $\pi : \tilde{X} \rightarrow X$ is chosen to be an equvariant desingularization, and $\tilde{\xi}$ the lift of $\xi$ to a holomorphic vector field on $\tilde{X}$.
 
\end{itemize}

To prove this one needs to note the following. Let $U$ be a neighbourhood of a point from $E$ over which $[\tilde{D}]$ is trivial. Let $d_U$ be the function which represents locally the metric $\tilde{d}$ there. Without loss of generality we may assume that $d_U$ is bounded in $U$. On $U\,-\,E$ the function $\varphi$ has the form $(1/k).\log(d_U.(t/\tilde{\omega})^k)$. If $\tilde{\xi}$ is as above, then the integral

$$\int_{U\,-\,E}\tilde{\xi}(\varphi)\,dV_{\tilde{\omega}} = \int_{U\,-\,E}\tilde{\xi}(\varphi).(\tilde{\omega}/t)\,dV_t$$ 
is small enough together with $U$.

\begin{theo}[Futaki\,-\,Ding\,-\,Tian, cf. \cite{D-T}, {\bf Lemma 1.2}] 
We have the following:

1) $F_{\varphi}$ does not depend on the choice of $\tilde{h}$, 

2) $F := F_{\varphi} : \eta (X) \rightarrow {\bf C}$ is a Lie algebras homomorphism.

\end{theo}

The morphism $F$ is called {\bf the Generalized Futaki invariant of} $X$. 
\par
We would like to emphasize that $F$ can be calculated by an integration over any complement to a rare subset of the smooth part $X_{reg}$ of $X$. This allows us to use two special representations of $F(\xi)$ which we are going to explain now.(cf., Mabuchi \cite{M1})
\par
Let $U \subset X$ be a Zariski\,-\,open subset, such that $\omega _X$ trivializes over $U\cap X_{reg}$, and let $b$ be the corresponding local frame for $\omega _X$. Then, the volume form $dV_{\omega}$ and the metric $l$ can be expressed as follows:

$$dV_{\omega} = \exp(-\,u_{\omega})(\lambda .b\wedge\bar{b}),\,\,\,\,l = \exp(-u_l)(\lambda.b\wedge\bar{b})^k,$$
where $\lambda := (\sqrt{-1})^n.(-1)^{n(n-1)/2}$. Hence, there is defined a volume form

$$dV_l := \exp(- u_l/k).\lambda.b\wedge\bar{b}.$$
In these notations we can express $\varphi$ over $U \cap X_{reg}$ as
$\varphi = u_{\omega} - (1/k).u_l$. Now, it is easy to calculate that

$$\xi (\varphi)\,dV_{\omega} = L_{\xi}(\varphi)\,dV_{\omega} = -\,L_{\xi}(dV_{\omega}) + div_l(\xi)\,dV_{\omega},$$
where $div_l(\xi)$ is the divergence of $\xi$ with respect to the volume form $dV_l$. The analytic Stokes' theorem (on $U$) gives the formula

$$\int_UL_{\xi}(\varphi)\,dV_{\omega} = \int_Udiv_l(\xi)\,dV_{\omega}.$$

By using this formula one can prove ( Mabuchi \cite{M2}, \S 6) that the nilpotent radical of $Lie(G)$ lies in $ker(F)$. ( In fact, in his proof Mabuchi uses admissible metrics, and so the same reason applies in the nonsmooth case.)  This means that if

$$Lie(G) = Lie(H) + Lie(Rad_uG)$$
is the Levi decomposition of $Lie(G)$, then $F$ is determined by its restriction to $Lie(H)$. The group H is a (maximal) reductive subgroup of G, and $F$ restricts to a character of $Lie(H)$. Hence, $F$ is determined completely by its restriction to $Lie(Z)$ - the Lie algebra of the centre of H.

To derive the second formula for $F$ suppose that $b\wedge\bar{b}$ is $\xi$-invariant. Then $div_l(\xi) = (-1/k).L_{\xi}(u_l)$ and we get the formula

$$\,\,\,\,\,\,\,\,\,\,\,\,\,\,\,F(\xi) = (-1/k.)\int_UL_{\xi}(u_l)\,dV_{\omega}.\,\,\,\,\,\,\,\,\,\,\,\,\,\,\,\,\,\,\,\,\,\,\,\,\,\,\,\,\,\,\,\,\,\,\, (\ast)$$
In particular, if $S$ is a subgroup of $G$, such that $U$ and $b\wedge\bar{b}$ are $S$-invariant, then the restriction of $F$ to the subalgebra $Lie(S)$ can be calculated by this formula (\cite{M1}).

The significance of the real part of the Generalized Futaki invariant stems from its connection with the existence of Einstein-K\"ahler metrics on Fano manifolds (via the notion of Weakly K-stability of Fano manifolds, due to G.Tian \cite{Tian}). We recall here some definitions, and state the main theorem concerning this connection.

 Let $\Delta$ be the unit disk in ${\bf C}^1$.\,\, The algebraic fibration without multiple fibres $\pi : W \rightarrow \Delta$\,\, is called a degeneration of the (Fano) manifold $M$,\,\, if for some $z\in \Delta$\, the fibre $W_z = \pi ^{-1}(z)$\,\, is biholomorphic to $M$.\,\, This degeneration is ${\bf special}$\, if its relative anticanonical sheaf is ample on the whole $W$,\,\,its central fibre $W_0 = \pi^{-1}(0)$ \,\, is a normal variety, and the  dilations $z \mapsto \lambda z$ \,\, on $\Delta\,\,(|\lambda| <1)$\,\, can be lifted to a family (1-parametric group) of automorphisms $\sigma (\lambda)$\,\, of $W$. \,\, In case $W \cong M\times\Delta$\,\, and $\pi$\,\, coincides with the projection onto $\Delta$,\,\, the degeneration is called {\bf trivial}.

For any special degeneration of a Fano manifold the central fibre $W_0$\,\, is an almost Fano variety. Moreover, the family of automorphisms $\sigma (\lambda)$\,\, determines a vector field $v_W := -\sigma'(1)$\,\,on $W$\,\, tangent to $W_0$ along its smooth part, hence an admissible vector field on $W_0$.

Denote by $F_{W_0}$\,\, the Generalized Futaki invariant of $W_0$.

\begin{defin} 
The Fano manifold $M$ is called {\bf Weakly K-stable} if for any special degeneration $W\,\, Re(F_{W_0}(v_W)) \geq 0 $,\,\, and equality holds if and only if $W$ \,\, is trivial.
\end{defin}

\begin{theo}[Tian \cite{Tian}]
If the Fano manifold $M$\,\,admits an Einstein-K\"ahler metric, then $M$\,\, is Weakly K-stable.
\end{theo}

\section{Moment maps - general setting.}

Let $(M, \omega)$ be a symplectic manifold; let $K$ be a compact group acting on $M$ by symplectomorphisms; let $k := Lie(K)$. The action of $K$ on $M$ gives a representation \,$R : k \rightarrow Lie(Aut(M))$\, as well as its dual \,$R^* : Lie(Aut(M))^* \rightarrow k^*$.

By definition (see \cite{Fut1}, Chapter VII), the map $\mu : M \rightarrow k^*$ is called a moment map of $M$ with respect to the action of $K$ if it has the following properties:

$$1) \forall g \in K\,\,\, \Rightarrow\,\,\, \mu\circ g = Ad(g^{-1})^*\circ \mu \,\,\,\,\,\,\,\,\,\,\,\,\,\,\,\,\,\,\,\,\,\,\,\,\,\,\,\,\,\,\,\,\,\,\,\,\,\,\,2)d\,\mu = R^*\circ \omega.\,\,\,\,\,\,\,\,\,\,\,$$

The classical theorem of Atiyah \cite{A1} explains the image of the orbits of $K$ under $\mu$ when $K$ is a d-dimensional torus : $K \cong S^1 \times \cdots S^1$ ($ d$ times). In the special case of a K\"ahler manifold $(M, \omega )$ and the torus $K$ acting by holomorphic isometries, it gives the following description :

\begin{itemize}
\item The complexification $G$ of $K$ acts holomorphically on $M$ too. Let $Y$ be an orbit of $G$ in $M$, and let $\overline{Y}$ be its closure. Hence, $\mu (\overline{Y})$ is a convex polytope $P$ in ${\bf R}^d$ ($\cong k^*$) and $\mu$ induces a homeomorphism of $\overline{Y}/K$ and $P$. This homeomorphism is a real analytic diffeomorphism of $Y/K$ and $Int(P)$ - the interior of $P$. If the complex dimension of $Y$ is $d$, then the mapping $Y \rightarrow Int(P)$ induced by $\mu$ is a principal $K$-bundle.

\end{itemize}

The most important example for us will be the case $M = {\bf P}^N$, $\omega$ - the Fubini\,-\,Studi form on it, and $G$ - the (N+1)\,-dimensional algebraic torus with the standard action on ${\bf P}^N$. In this case the compact group is (N+1)-\,dimensional compact torus $K = S^1\times\cdots\times S^1$, the corresponding Lie algebra is generated by (N+1) elements : $k = e_1.{\bf R}\,+\,\cdots\,+\,e_N.{\bf R}$, and the j-th coordinate of the moment map , up to an additive constant, is given by the formula $\mu_j = {\mid x_j\mid ^2 \over \sum_{k=0}^{N}\mid x_k\mid ^2}$. In other words, if $\{e_0^*,\ldots,e_N^*\}$ is the basis of $k^*$ dual to the chosen basis of $k$, then

$$\mu (x_0:\ldots:x_N) = \sum_{i=0}^{N}{\mid x_i\mid^2\over\sum_{k=0}^{N}\mid x_k\mid^2}.e_i^*.$$

We will use this concrete moment map to induce moment maps on almost Fano varieties in the following way (for details see \cite{Ful}).

Let $X$ be a variety on which the $d$\,-dimensional algebraic torus $H = {\bf C}^*\times\cdots{\bf C}^*$ acts holomorphically. Let $L$ be a very ample line bundle on $X$, such that the acton of $H$ on $X$ lifts to a linear action on $L$. Hence $H$ acts on the space $\Gamma(X, L)$ of global sections of $L$ and we can find a basis of semi-invariants of it

$$ \{\{s_0,\cdots, s_N\}\,\,|\,\, h(s_i) = \exp(\lambda_i(h)).s_i, \,\,\,i = 0,\dots, N\}.$$

Let $\Phi : X \hookrightarrow {\bf P}^N$ be the embedding of $X$ determined by $L$. In the basis ${s_0,\ldots, s_N}$ this embedding defines a representation $R : H \rightarrow G$ together with the Lie algebra representation $dR : Lie(H) \rightarrow Lie(G)$ and its dual $dR^* : Lie(G)^* \rightarrow Lie(H)^*$.  Now, if we choose the preimage of $K$ under $R$ to be the maximal compact subgroup $T$ of $H$, then $dR^*$ restricts to a representation $\rho : k^* \rightarrow t^*$. The composition map $\mu_X := \rho \circ \mu \circ \Phi$\,\, is in fact a moment map of $X$ with respect to the action of $T$ on it. We have :

$$\mu_X(x) = \sum_{i=0}^{N}{\mid s_i(x)\mid^2 \over \sum_{k=0}^{N}\mid s_k(x)\mid^2}.\rho(e_i^*).$$

The Atiyah theorem has a very nice description when $X$ is a toric variety. We will explain it in the next section.

\section{The Generalized Futaki invariant of almost Fano toric varieties.}

Let $X$ be a normal complex variety of dimension $n$. Let $T^n$ be an $n$\,-dimensional algebraic torus.

\begin{defin} $X$ is called a toric variety if there exist an embedding $\varphi : T \hookrightarrow X$ and an action $R : T\times X \rightarrow X$ such that $\varphi (T)$ is an open dense subset, and for any $t_1, t_2 \in T$ we have $R(t_1, \varphi(t_2)) = \varphi (t_1.t_2)$.
\end{defin}

Since ${\bf C}^*\times\ldots\times{\bf C}^* = Hom_{{\bf Z}}(N, {\bf Z})$ for some lattice $N \cong {\bf Z}^{\oplus n}$, then to choose an isomorphism $T \cong {\bf C}^*\times\ldots\times{\bf C}^*$ is equivalent to choosing a lattice $N$ and an identification $T = Hom_{\bf Z}(N, {\bf C}^*)$.
\par
{\bf Fact.} Having chosen an identification $T = Hom_{\bf Z}(N, {\bf C}^*)$, there exists a unique {\bf fan} $\Delta$ in $N_{\bf R} := N\otimes_{{\bf Z}}{\bf R}$ which determines $X$.
\par
 
$\Delta \subset N_{\bf R}$ is called a fan if
\begin{itemize}
\item $\Delta$ is a collection of cones in $N_{\bf R}$ with vertices at the origin of $N_{\bf R}$;

\item each cone of $\Delta$ is generated by a finite number of vectors from $N$, and does not contain lines through the origin;

\item the faces of the cones of $\Delta$ belong to $\Delta$;
\item the intersection of any two cones of $\Delta$ is a face of each.
\end{itemize}

The union, $\mid\Delta\mid$, of all cones of $\Delta$ is called the {\bf support} of $\Delta$.
\par
The set of all $n-$dimensional cones of $\Delta$ \,\, is denoted by $\Delta^n$.
\par
Any cone $\sigma \in \Delta$ \,\, is a convex hull of its one-dimensional extremal rays. If $\tau$ \,\, is a ray, then $\tau \cap N$ has a unique generator $v$ \,\,over ${\bf N}$.\,This generator is called {\bf the primitive generator} of $\tau$.
\par 

The reconstruction of $X$ by the fan $\Delta$ goes as follows.
Let $ M := Hom_{{\bf Z}}(N, {\bf Z})$ be the dual lattice of $N$, and let $<\,,\,> : M\times N \rightarrow {\bf Z}$ be the corresponding pairing between them. For any cone $\sigma \in \Delta$ there is defined the {\bf dual cone} $\check{\sigma} \subset M_{\bf R}$ as follows:

$$\check{\sigma} := \{ u \in M_{{\bf R}} : <u, v> \geq 0 \,\,\forall v \in \sigma\}.$$
The set of elements $S_{\sigma} := \check{\sigma} \cap M$ form a semi-group; the corresponding group algebra ${\bf C}[S_{\sigma}]$ is finitely generated and hence defines an affine scheme $U_{\sigma} := Spec({\bf C}[S_{\sigma}])$ which is called {\bf the affine toric variety, corresponding to} $\sigma$. If $\tau := \sigma \cap \sigma_1$ is a face of two cones of the fan, then in a natural way $U_{\tau} \subset U_{\sigma}$ and $U_{\tau} \subset U_{\sigma_1}$. Hence, we can glue $U_{\sigma}$ and $U_{\sigma_1}$ along $U_{\tau}$. Glueing together the elements of $\{U_{\sigma} : \sigma \in \Delta\}$ we reconstruct $X$. In particular, the torus $T \hookrightarrow X$ corresponds to the cone ${0} \in \Delta$, and its dual cone coincides with the whole $M_{{\bf R}}$. Hence,

$${\bf C}[S_{{0}}] = {\bf C}[t_1,t_1^{-1},\ldots,t_n,t_n^{-1}],$$
where $t_i$ corresponds to the basis element $e_i^*$ ($i = 1,\ldots,n$). 
\par
Each element $u$ of the lattice $M$ determines a holomorphic function $\chi^u$ on $T$: if $ u= \sum_{i=1}^{n}\alpha_i.e_i^*,\,\,\alpha_i \in {\bf Z}, \,\,i=1,\ldots,n$,\,\,\,then $\chi^u := t_1^{\alpha_1}.\cdots.t_n^{\alpha_n}$. The functions $\chi ^u$ can be extended to rational functions on $X$.
\par
The fan $\Delta$ carries the information about the properties of $X$. For example
\begin{enumerate}
\item $X$ is compact iff $\mid\Delta\mid = N_{\bf R}$;
\item $X$ is smooth iff each cone in $\Delta$ is generated by a subsystem of generators of $N$;
\item the (n-k)\,-dimensional irreducible $T$-invariant subvarieties of $X$ correspond bijectively to the k\,-dimensional cones of $\Delta$. In particular, the $T$\,-invariant irreducible effective Weil divisors of $X$ correspond bijectively to the rays of $\Delta$ : $\tau_1,\ldots,\tau_d$\,\, (with corresponding primitive generators over $v_1,\ldots,v_d$)
\end{enumerate}

The {\bf dualizing sheaf} of $X$ is defined to be the sheaf $\omega_X := {\cal O}_X(- \sum_{j=1}^{d}D_j)$, where $D_j$ is the Weil divisor, corresponding to $\tau_j,\,\, j = 1,\ldots ,d$.

If $D$ is a Cartier divisor on $X$, then, since $X$ is normal, it defines a Weil divisor $[D]$. If $D$ is a $T$\,-Cartier divisor, then $[D] = \sum_{i=1}^{d}a_i.D_i \,,\, a_i \in {\bf Z}, i=1,\ldots,d$. Such a $T$\,-Cartier divisor determines a rational convex polytope (with vertices in $M$):

$$P_D := \{ u \in M_{{\bf R}} : <u, v_i> \geq - a_i \,\, i = 1,\ldots,d \}.$$

{\bf Fact.} For each $T$\,-Cartier divisor $D$

$$\Gamma(X, {\cal O}_X(D)) = \oplus_{u\in P_D\cap M}{\bf C}.\chi^u.$$

We are primarily interested in the case of almost Fano toric varieties $X$. This means that $X$ is compact, and for some positive integer $k$\, the sheaf
$${\cal O}_X(k.\sum_{j=1}^dD_j)$$
is an ample invertible sheaf on $X$. In particular, since the toric varieties are Cohen-Macaulay,  $X$ is a {\bf Q}\,-Gorenstein variety.

The property of being {\bf Q}\,-Gorenstein toric variety in terms of the properties of the corresponding fan is interpreted as follows:
\begin{itemize}
\item For each n\,-dimensional cone $\sigma \in \Delta^n$ there exists an element $k_{\sigma} \in M_{\bf Q}$, such that $<n_i, k_{\sigma}> = 1$ for all primitive generators $n_1,\ldots, n_r$ of $\sigma$.
\end{itemize}

$X$ is Gorenstein iff $k_{\sigma} \in M$ for all n\,-dimensional cones $\sigma \in \Delta^n$.

\begin{itemize}
\item A {\bf Q}\,-Gorenstein toric variety $X$ is almost Fano iff the convex hull
$$ P_{-K_X} := coh< - k_{\sigma} : \sigma \in \Delta^n> $$
is an n\,-dimensional polytope with vertices exactly $k_{\sigma}$ \,\, for all $\sigma \in \Delta^n $.
\end{itemize}
In this case the polytope corresponding to the line bundle $L$ is $k.P_{-K_X}$.

Suppose $L = \sum_{i=1}^d k.D_i$ is a very ample T-Cartier divisor on the almost Fano toric variety $X$. Let ${u_o = o, u_1,\ldots, u_N}$ be the set of integer points of the corresponding polytope $k. P_{-\omega_X} \subset M_{\bf R}^X$. Then there is an embedding $\phi : X \hookrightarrow {\bf P}^N$ given by the mapping

$$ X \ni x \longmapsto (\chi^{u_0}(x):\dots:\chi^{u_N}(x)).$$

If ${b_0^*, \ldots, b_N^*}$ is a basis of $Lie(S^1\times\dots\times S^1)^* = M_{\bf R}^{{\bf P}^N}$, then the moment map

$$\mu^{{\bf P}N} : {\bf P}^N \longrightarrow M_{\bf R}^{{\bf P}N}$$
is given by the formula

$$ \mu^{{\bf P}^N}(x_0:\dots:x_N) = {1\over \sum_{i=0}^N|x_i|^2}.\sum_{i=0}^N|x_j|^2.b_j^*$$

The maximal compact subtorus $C_T := S^1\times\dots\times S^1 \subset T$ ( $n$ times) is given by $\{e^{i\theta_1}\times\dots\times e^{i\theta_n} | \theta_i \in (-\pi, \pi)\}$. Accordingly, the restriction of the morphism $\phi$ to this torus is given by the formula

$$e^{i\theta_1}\times\dots\times e^{i\theta_n} \longmapsto (e^{i\sum\theta_j\alpha_j(u_0)} :\dots : e^{i\sum\theta_j\alpha_j(u_N)}).$$
Hence, $\phi$ induces a morphism of tori

$$\varphi : C_T \longrightarrow C := S^1\times\dots\times S^1 \subset {\bf C}^*\times\dots\times{\bf C}^* \,\, (N+1 times)$$
with corresponding algebra morphism $d\varphi : Lie(C_T) \rightarrow Lie(C)$ given by $d\varphi (e_i) = \sum_{j=0}^N \alpha_i(u_j).b_j$.

We will consider the moment map of $X$ induced by $\mu^{{\bf P}^N}$ : $$\mu^X := d\varphi^*\circ\mu^{{\bf P}^N}\circ\phi.$$ It is easily checked that the k-th coordinate of this map (in the basis ${e_1^*,\ldots ,e_n^*}$) is given by the formula

$$\mu^X_k (x) = {\sum_{j=0}^N\alpha_k(u_j).|\chi^{u_j}(x)|^2\over\sum_{i=0}^N|\chi^{u_i}(x)|^2} , \,\,\, k = 1,\ldots, n.$$

Let us restrict $\mu^X$ to the torus $T \subset X$, and let introduce the polar coordinates for $t_j \in {\bf C}^*$ :

$$ t_j := \exp(x_j +\sqrt{-1}\theta_j),\,\, x_j \in {\bf R}, \theta_j \in (-\pi, \pi), j=1,\ldots,n.$$
It is easily verified that

$$\mu^X_k (t_1\times\dots\times t_n) = (1/2).\partial_k f(x_1,\ldots,x_n),$$
where $$f(x_1,\ldots,x_n) := \log (\sum_{i=0}^N|\chi^{u_i}(t_1\times\dots\times t_n)|^2),$$ and $\partial_k$ denotes the partial derivative with respect to $x_k$.
\par
Returning to the embedding $\phi : X \rightarrow {\bf P}^N$ we notice that the k-th tensor power of the  T-invariant n-form $b := {dt_1\over t_1}\wedge\dots\wedge {dt_n\over t_n}$ on $T$ corresponds to $ 0 \in k.P_{-\omega_X}$, and in the local frame $b^{\otimes k}$ of the very ample line bundle $L$ \,\,the map \,\,$\phi$ is given by $$X_0 =1,\,\, X_j = \prod_{i=1}^n t_i^{\alpha_i(u_j)},\,\, j=1,\ldots,N.$$
\par
Choose the metric $\tilde{h} = (\sum_{i=0}^N|X_i|^2)^{-1}$ on the hyperplane line bundle $[H]$ on ${\bf P}^N$. Hence locally in the frame $b^{\otimes k}$ the metric $l$ on $L$ is given by the function $(1 + \sum_{i=1}^N|\prod_{j=1}^n|t_j^{\alpha_j(u_i)}|^2)^{-1}$, and we get that the corresponding function $u_l$ is equal to $f$.

Now we are ready to obtain the relation between the generalized Futaki invariant and the moment map of $X$.\, First of all ${\partial\over\partial x_k}$ identifies with the vector field $t_k{\partial\over\partial t_k} + \bar {t}_k{\partial\over\partial\bar{t}_k}$, and hence

$$ F({\partial\over\partial x_k}) = 2.Re (F(t_k{\partial\over\partial t_k})).$$

Since  $b\wedge\bar{b}$ is T-invariant we can use the formula (*) and to see that 
$$F({\partial\over\partial x_k}) = -(1/k).\int_{T}{\partial u_l\over\partial x_i}\,dV_{\omega} = -(2/k).\int_{T}\mu_i\,dV_{\omega}.$$

On the other hand, $\omega = (1/k).\phi^*\omega'$, and $\omega' = Ric(\tilde{h})$, which gives that

$$\omega = (-1/k).\sqrt{-1}.\partial\bar{\partial}f = (1/4k)\sqrt{-1}\sum_{i,j}{\partial^2f\over\partial x_i\partial x_j}\,{dt_i\over t_i}\wedge{d\bar{t}_j\over \bar{t}_j}.$$
In polar coordinates introduced above we can express the volume form as follows:

$$dV_{\omega} = (1/n!).\omega^{\wedge n} = (1/4k)^n\det{\partial^2f\over\partial x_i\partial x_j}\prod_{i=1}^n(\sqrt{-1}{dt_i\over t_i}\wedge{d\bar{t}_i\over \bar{t}_i}) $$

$$= (1/2k)^n\det{\partial^2f\over\partial x_i\partial x_j}\prod_{i=1}^ndx_i\wedge d\theta_i.$$

We can conclude now the calculation of the generalized Futaki invariant:

$$(-k).F({\partial\over\partial x_s}) = \int_{T}{\partial f\over\partial x_s}\,dV_{\omega} = (1/2k)^n\int_{T}{\partial f\over\partial x_s}\,\det{\partial^2f\over\partial x_i\partial x_j}\prod_{i=1}^ndx_i\wedge d\theta_i $$

$$= (\pi/k)^n\int_{{\bf R}^n}{\partial f\over\partial x_s}\det{\partial^2f\over\partial x_i\partial x_j}\,dx_1\wedge\ldots\wedge dx_n $$

$$= (\pi/k)^n\int_{{\bf R}^n}{\partial f\over\partial x_s}\,d{\partial f\over\partial x_1}\wedge\ldots\wedge{\partial f\over\partial x_n} $$

$$= (\pi/k)^n\int_{2\mu^X({\bf R}^n)}y_s\,dy_1\wedge\ldots\wedge y_n.$$

Finally we get the formula:

$$Re F(t_s{\partial\over\partial t_s}) = -(2\pi)^n\int_{P_{-\omega_X}}y_s\,dy_1\wedge\ldots\wedge dy_n,$$
where $y_1,\ldots, y_n$ are the coordinates of $M_{\bf R}^X$\, in the basis $e_1^*,\ldots, e_n^*$..

This last formula allows us to find explicitly the real part of the restriction to the torus $T$ of the Generalized Futaki invariant of $X$.

\begin{rem}
The numbers $F({\partial\over\partial x_s}),\,s=1,\ldots, n$ are in close relation with the barycentre of the polytope $P_{-\omega_X} \subset M_{{\bf R}}^X$ with respect to the natural Euclidean measure in $M_{{\bf R}}^X$. More precisely, denote by ${\bf a} = (a_1,\ldots, a_n)$ the baricentre of $P_{-\omega_X}$\,, i.e., the element of $M_{{\bf R}}^X$\, with coordinates

$$a_s = {\int_{P_{-\omega_X}}y_s\, dy_1\wedge\ldots\wedge dy_n \over \int_{P_{-\omega_X}}dy_1\wedge\ldots\wedge dy_n}.$$
The calculations above give that

$${1\over 2}.F({\partial\over\partial x_s}) = - (2\pi.c_1(X))^n.a_s,\,\,\,\,s=1,\ldots, n.$$
Note that $-{\bf a}$ is the barycentre of $-P_{-\omega_X}$.
\end{rem}

\begin{rem}
It is well known that if $(X, T)$ \,is a smooth toric variety, then $T$\, is a maximal torus of $Aut(X)$,\,( cf. \cite{Oda}, \S 3.4). For a nonsmooth toric variety $(X, T)$ \,let $\varphi : \tilde{X} \rightarrow X$\, be its canonical resolution of singularities. Then, we have a natural inclusion $Aut(X) \hookrightarrow Aut(\tilde{X})$\,, and that, identifying $T$\, with its image under this inclusion, $(\tilde{X}, T)$ is a toric manifold. From this it follows that $T$\, is a maximal torus of both groups - $Aut(\tilde{X})$\, and $Aut(X)$. Hence, in the Levi decomposition of $Aut(X)$\, we can choose a maximal reductive subgroup which has $T$\, as a maximal torus. 
\end{rem}

We can summarize the above considerations in the following Theorem

\begin{theo}\label{torus}

Let $X$ be an almost Fano variety which is toric with respect to the action of a torus $T \subset Aut(X) =: G$.\,\, Fix an isomorphism $T \cong Hom(N, {\bf C}^*)$,\,\, where $N$ is a lattice with dual lattice $M$,\, and let $P_{-\omega_X}$ \, be the polytope in $M_{{\bf R}}$\,\, corresponding to the dualizing sheaf $\omega_X$\,\, of $X$.\, Let ${\bf a}$\, be the barycentre of $P_{-\omega_X}$ as above. Suppose  $\xi \in Lie(G),\,\, \eta \in Lie(T) \subset Lie(G)$\,\, are such that

$$\xi - \eta \in Lie(Rad_uG) + Lie([H, H]) \, \subset \, Lie(G),$$
where $H$ is a maximal reductive subgroup of $G$ for which $T$\, is a maximal torus.

Then, $F(\xi) = F(\eta)$,\,\,where $F$ \, is the Generalized Futaki invariant of $X$.\,\, Moreover, if $\eta = \sum_{s=1}^n\eta^s.t_s(\partial/\partial t_s)$ \,\,in $Lie(T)$,\,\,then

$$\,\,\,\,\,\,\,\,\,\,\,\,\,\,\,Re(F(\xi)) = - (2\pi.c_1(X))^n. \sum_{s=1}^n Re(\eta^s).a_s.\,\,\,\,\,\,\,\,\,\,\,\,\,\,\,\,\,\,\,\,\,\,\, (**)$$

\end{theo} 

Indeed, in polar coordinates $t_j = \exp(x_j + \sqrt{-1}\theta_j)$\, we have
$t_j\partial/\partial t_j = (1/2).(\partial/\partial x_j - \sqrt{-1}.\partial/\partial \theta_j)$\,, and hence
$$F(\sum_{s=1}^n\eta^s.t_s\partial/\partial t_s) = {1\over2}.\eta^s.(F(\partial/\partial x_s) - \sqrt{-1}.F(\partial/\partial \theta_s)).$$
Since $\tilde{h}$\, in the definition of the Generalized Futaki Invariant of $X$\, is $C_T$-\,invariant, then $F(\partial/\partial \theta_s) = 0$\,. This proves $(**)$ from the Theorem above.

\section{Examples.}

\subsection{Preliminaries}

Our examples concern the degenerations of the three-dimensional Fano manifold $V_{38}$\,, the blowup of the projective space in a twisted cubic curve, into irreducible normal toric varieties. For each twisted cubic curve $C \subset {\bf P}^3$\,\, the corresponding blowup $B_C{\bf P}^3$\,\,can be represented as a smooth intersection of the Segre embedding of ${\bf P}^3\times{\bf P}^2$ \,(in ${\bf P}^{11}$)\, by a codimension-two hyperplane. Explicitly this can be shown as follows.

Let ${\bf P}^3$\, be the projectivization of a 4-dimensional vector space $V$,\, and ${\bf P}^2$\, be the projectivization of a 3-dimensional vector space $W$. Then, ${\bf P}^{11}$ can be identified with the projectivization of the tensor product of $V$ \, and $W$ \,: ${\bf P}^{11} \cong P(V\otimes W)$.

Fix coordinates $(X_0:X_1:X_2:X_3)$ \, in ${\bf P}^3$. Then, the ideal $I$\, of the twisted cubic curve $C$ is generated by $(2\times 2)$\,-minors of \, \,$(2\times 3)$\, matrix $A$\, with entries - linear forms on ${\bf P}^3$\,:

$$A = \left( \begin{array}{ccc}
a_{10} & a_{11} & a_{12} \\
a_{20} & a_{21} & a_{22}
\end{array} \right),$$

$$I=(\Delta_1, \Delta_2, \Delta_3),$$
where $a_{ij}$\, are linear forms on ${\bf P}^3$,\,\, and $\Delta_i$\,\, is the $(2\times2)$-minor of $A$\, that does not contain its  $i$th column. Moreover, for each pair of complex numbers $(\lambda,\mu) \not= (0,0)$\,\, the linear forms

$$\lambda.a_{10} + \mu.a_{20},\ \ \lambda.a_{11} + \mu.a_{21},\ \ \lambda.a_{12} + \mu.a_{22}$$
are linearly independent. 

If we denote by $H_i$\, the hyperplane of ${\bf P}^{11}$\, defined by the equation $$h_i({\bf X},{\bf Y}) = a_{i0}({\bf X}).Y_0 + a_{i1}({\bf X}).Y_1 + a_{i2}({\bf X}).Y_2,$$ 
then $H_1$ \,and $H_2$\, are linearly independent and

$$ B_C{\bf P}^3 \cong ({\bf P}^3\times{\bf P}^2) \cap H_1 \cap H_2.$$

Our purpose here is to investigate the degenerations of $V_{38}$\, obtained by moving the codimension-two hyperplane and intersecting with ${\bf P}^3\times{\bf P}^2$. The reason for doing so is that for any degeneration $\pi : W \rightarrow \Delta$\, of $V_{38}$\,in ${\bf P}^3\times{\bf P}^2$\, its relative anticanonical sheaf is $\pi$-very ample, so if $\sigma(\lambda)$\, is an one parametric group as required in the definition of special degenerations, and $W_0$\, is normal, then $W$\, is automatically a special degeneration of $V_{38}$.

\subsubsection{}

Let ${\bf P}^9$ be a subspace of ${\bf P}^{11}$. This subspace is the zero set of two linear forms on ${\bf P}^{11}$,\, which can be represented by a matrix $A$\, as above. This matrix defines an ideal $I$\, generated by its maximal minors. 

The intersection ${\bf P}^9 \cap ({\bf P}^3\times {\bf P}^2)$\, has dimension 3 or 4. One of its irreducible components is isomorphic to the blowup of ${\bf P}^3$\, in the ideal $I$. And since we are interested in irreducible intersections only, we have to avoid considering ${\bf P}^9$'s defined by linear forms of the form
$$h_i = f({\bf X}).h_i'({\bf Y}),$$
where $h_1',\, h_2'$\, are linearly independent forms on $W$,\, and $f$ \, is a nonzero linear form on $V$. These last form a 5-dimensional closed subset (isomorphic to ${\bf P}^3\times{\bf P}^2$) of the Grassmannian $Gr({\bf P}^9, {\bf P}^{11})$.

On the other hand, the three-dimensional intersections ${\bf P}^9 \cap ({\bf P}^3\times{\bf P}^2) \subset {\bf P}^{11}$,\, are characterized as those subschemes $Y \subset {\bf P}^3\times {\bf P}^2$\,\, which have Hilbert polynomial

$$P(n) := \chi ({\cal O}_Y(n,n)) = {1\over 3!}.(10n^3 + 24n^2 + 20n + 6),$$
and satisfy the vanishing condition

$$h^q({\cal I}_Y(1,1)) = 0, \ \ \ q\geq 1.$$

The last condition implies that the three-dimensional intersections with ${\bf P}^9$ form an open subset of the Hilbert scheme $Hilb^{P(n)}({\bf P}^3\times{\bf P}^2)$.  In other words, we are interested in the irreducible normal degenerations of $V_{38}$,\, which lie in the closure of the subset of the Hilbert scheme above corresponding to the smooth intersections ${\bf P}^9 \,\cap({\bf P}^3\times{\bf P}^2)$.

\subsubsection{}

Denote by $p : {\bf P}^9\,\cap({\bf P}^3\times{\bf P}^2) \,\rightarrow {\bf P}^3$\,\, the restriction of the natural projection ${\bf P}^3\times{\bf P}^2 \rightarrow {\bf P}^3$.\,Then, $p$ is biregular outside the preimage of the subscheme on ${\bf P}^3$\, defined by $I$,\, while its restriction to this preimage has positive dimensional fibres (which are linear subspaces of ${\bf P}^2$).\, Hence, if the intersection above is irreducible, then $I$ \, defines one-dimensional subscheme of ${\bf P}^3$.

Let $S$\, be a two-dimensional vector space. Denote by $E$ \, the set of $(2\times 3)$ matrices with entries - linear forms on $V$ \,:\, $E = Hom_{\bf C}(W\otimes V^*, F)$.\,Let, as before, $I = (\Delta_1, \Delta_2, \Delta_3)$ \, denote the ideal generated by the $(2\times2)$- minors of the elements of $E$.

 The group $G_1 := GL(W)\times GL(F)$\, acts on $E$ via
$$(g,h)A := hA(g\otimes id_{V^*})^{-1}.$$
This action obviously descends to an action of $G := G_1 /\Gamma $,\, where 
$$\Gamma := \{(\alpha.id_{W},\alpha.id_{F}) | \alpha \in {\bf C}^*\}.$$
\par$ $\par
The proofs of the following facts can be found in \cite{E-P-S}:
\begin{enumerate}
\item If $I$\, defines a curve in ${\bf P}^3$,\, then $\dim_{{\bf C}} <\Delta_1, \Delta_2, \Delta_3> = 3$.
\item Denote by $X$ the nets of quadrics in ${\bf P}^3$\, generated by $(2\times 2)$- minors of a matrix from $E$,\, and by $U \subset E$\, the set of all matrices for which $\dim_{{\bf C}} <\Delta_1, \Delta_2, \Delta_3> = 3.$ \, Then $X$ \,is smooth, and $U$ \,has a structure of a principal $G$-bundle over $X$.
\item Let $H$ be the component of $Hilb^{3t+1}({\bf P}^3)$\,\, containing the points corresponding to the twisted cubic curves, and let $U_1 \subset U$\,consist of matrices which ideals $I$ \, do not define curves. Then, $U_1$ maps onto  the subset $X_1$ of $X$ \, of ``point-plane incidence correspondence'', and there is a morphism $f : H \rightarrow X$\, which is the blowup of $X$ along $X_1$. 

\end{enumerate}

On the other hand, for each ${\bf P}^9 \subset {\bf P}^{11}$\, different choices of hyperplanes defining it determine matrices which differ by the action of (the image of) $GL(F)\times\{id_W\}$\, (in $G_1$).

Let $\tilde{U} \subset Gr({\bf P}^9, {\bf P}^{11})$\, denote the set of ${\bf P}^9$'s defined by elements of $U$. Then, $\tilde{U}$\, has a structure of a principal $PGL(3)$- bundle over $X$. Denote by $\phi : \tilde{U} \rightarrow X$\, the corresponding morphism. This morphism defines a morphism $\psi : \tilde{U} - \phi^{-1}(X_1) \rightarrow H - f^{-1}(X_1)$.\,\,Hence,  $\tilde{U}_1 := \tilde{U} - \phi^{-1}(X_1)$\,\, is the set of those ${\bf P}^9$'s which ideals $I$\, define curves.

The set $\tilde{U}$\, can be naturally identified with an open subset $U'$\,of\\ 
$Hilb^{P(n)} ({\bf P}^3 \times {\bf P}^2)$.\, Let $U''$ \,denote the subset corresponding to $\tilde{U}_1$\, in this identification.

\begin{prop}\label{mor}
There is a morphism $\psi : U'' \rightarrow H$\,  which is a principal $PGL(3)$- bundle over its image $H''$.
\end{prop}

 Note that the action of $PGL(3)$\, on $U''$\, corresponds to the automorphisms of ${\bf P}^{11}$\, which leave ${\bf P}^3\times{\bf P}^2$\, invariant.

\subsubsection{}

Now we are in a position to explain the relation between flat degenerations of $V_{38} \subset {\bf P}^3\times {\bf P}^2 \subset {\bf P}^{11}$ \,we are interested in, and the degenerations of the corresponding curves in ${\bf P}^3$.

Suppose we are given a flat family $\tilde {Y} \subset \Delta \times ({\bf P}^3\times {\bf P}^2)$\,\, in ${\mathit Hilb}_{{\bf P}^3\times{\bf P}^2}^{P(n)}(\Delta)$,\, which corresponds to a morphism $y : \Delta : \rightarrow U''$.\,\, Then the composition
$$\psi \circ y : \Delta \rightarrow H''$$
gives in turn a flat family $\tilde {V} \subset \Delta\times{\bf P}^3$\,\, in ${\it Hilb}_{{\bf P}^3}^{3n+1}(\Delta)$.

Consider the natural projection

$$\pi : \Delta\times({\bf P}^3\times{\bf P}^2) \rightarrow \Delta\times{\bf P}^3.$$
If $\tilde {Y}_{\lambda}$\, is irreducible ($\lambda \in \Delta$)\,, then the restriction of $\pi$\, to $\tilde {Y}_{\lambda}$\, is the blowup of ${\bf P}^3$\, in $\tilde {V}_{\lambda}$. Hence, if for all $\lambda \in \Delta$\, the fibre\,\,$\tilde {Y}_{\lambda}$\, is irreducible, then the restriction of $\pi$ \, to $\tilde {Y}$\,\, is just the blowup of $\Delta\times{\bf P}^3$\,\, in $\tilde {V}$.

We are primarily interested in flat families $\tilde {Y}$\,which are degenerations of $V_{38}$\, with an action of 1-parametric group of transformations $\{g_t\}_{t\in (0,1]}$\,\, of $\tilde {Y}$\, which lifts the dilation automorphisms ,\,\,$ \lambda \mapsto t.\lambda$\,\,, of \, $\Delta$:

$$g_t|_{\tilde {Y}_{\lambda}} : \tilde {Y}_{\lambda} \rightarrow \tilde {Y}_{t\lambda}.$$
This group of transformations leaves the central fibre $\tilde {Y}_0$\, invariant.

When $\tilde {Y}_{\lambda}$\,\, is irreducible for all $\lambda \in \Delta$\,\,the group\, $\{g_t\}$\,  descends (under $\pi$\,) to a group of transformations $\{h_t\}$\,\,of $\Delta \times{\bf P}^3$\,\, which preserves $\tilde {V}$\,(lifting the dilation automorphisms of $\Delta$\,). Hence, $\{h_t\}$\, is a 1-parametric group in
$$Aut(\tilde {V}_0 \subset {\bf P}^3) = Aut(\tilde {Y}_0).$$
In this case $\tilde {V}_{\lambda},\,\,(\lambda \not= 0)$\,\, is a twisted cubic curve in ${\bf P}^3$.

Conversally, suppose we are given a flat family $\tilde {V}_{\lambda} \subset \Delta\times{\bf P}^3$\, in ${\it Hilb}_{{\bf P}^3}^{3n+1}(\Delta)$\,\, with an action of 1-parametric group $\{h_t\} \subset Aut({\bf P}^3)$\,\, such that
$$h_t |_{\tilde {V}_{\lambda}} : \tilde {V}_{\lambda} \rightarrow \tilde {V}_{t\lambda},$$
all $\tilde {V}_{\lambda} \,\,(\lambda \not=0)$\, are twisted cubic curves, and $\tilde {V}_0$ \, corresponds to irreducible intersection of ${\bf P}^3\times{\bf P}^2$\, with ${\bf P}^9$.\,\,Denote by $z : \Delta \rightarrow H''$\,\, the morphism to which the family $\tilde {V}$\, corresponds. Such a flat family defines a degeneration of $V_{38}$\, together with a group action, which lifts the dilation morphisms of $\Delta$\,, in the following way.

Let $\tilde {Y}'$\, be the blowup of $\Delta\times{\bf P}^3$\, in $\tilde {V}$. \,Hence, $\tilde {Y}'$ is a flat family, with fibres $\tilde {Y}'_{\lambda}$ isomorphic to the blowup of ${\bf P}^3$\, in $\tilde {V}_{\lambda}$,\, to which the action of $\{h_t\}$\, lifts (defining a group $\{g_t\}$)\,. As a flat family of $V_{38}$'s\, with a group  acting on it, the pair \,\,$(\{\tilde{Y}'_{\lambda \not= 0}\}, g_t)$\,\, can be represented, up to an action of $ PGL(3)$, as  intersection of \,\,${\bf P}^3\times{\bf P}^2$\,\, with a family, $\{H_t\}$\,\, of \,${\bf P}^9$'s\, in \,${\bf P}^{11}$\,, over \,$\Delta$\,,\,with a group \,$\{\tilde{h}_t\} \subset Aut({\bf P}^{11})$\,\, acting on it, lifting the action of \,$\{h_t\}$\,\,(under the morphism \,$\psi$\,\,from Proposition \ref{mor}).\, \,\,Indeed, fix \,$\lambda_0 \in \Delta^*$\, and let \,$A_{\lambda_0} \in \psi^{-1}(z(\lambda_0))$\,\, be a \,$(2\times3)$\,\, matrix, with entries - linear functions on \,${\bf P}^3$\,\,, which maximal minors determine the ideal of \,$\tilde{V}_{\lambda_0}$\,\,as a subscheme of \,${\bf P}^3$\,. Such a matrix is unique up to the action of the group \,$G$\, (from  Subsubsection 5.1.2). The group \,$\{h_t\}$\, defines then a family \,$\{A_{\lambda}\}$\,, on which it acts lifting the dilation morphisms of \,$\Delta^*$.\, Let \,$A_0$\, be the limit \,\,$\lim_{t \rightarrow 0} h_t(A_{\lambda_0})$\,\,. Note that, for \,$\lambda \not= 0$\,, the maximal minors of \,$A_{\lambda}$\, define the ideal of the subscheme \,\,$\tilde{V}_{\lambda} \subset {\bf P}^3$\,.\, It follows from this that, if \,$H_{\lambda}$\, denotes the codimension-two hyperplane of \,${\bf P}^{11}$\, defined by \,$A_{\lambda}$\,, then the intersection \,\,$({\bf P}^3\times{\bf P}^2) \cap H_{\lambda}$\,\, is isomorphic to \,$\tilde{Y}_{\lambda}$\, for \,$\lambda \not= 0$\,.\, Furthermore, the natural representation of the group \,$\{h_t\}$\, into \,$Aut({\bf P}^{11})$\,\,\,(under the identification \,\,${\bf P}^{11} = {\bf P}(V\otimes W)$\,\,) defines a group \,$\{\tilde{h}_t\}$\, which acts on \,$\{H_{\lambda}\}$\, lifting the dilation morphisms of \,$\Delta$\,, as was asserted. The family \,$H_t$\, depends only on the choice of \,$H_{\lambda_0}$\,, and hence is unique up to action of \,$ PGL(3)$. The restriction of \,$\tilde{h}_t$\,\, to\, $\tilde{Y}'_{\lambda}$\,\, coincides with \,$g_t$.

Denote by \,$\tilde{Y}$\, the degeneration of \,$V_{38}$\, corresponding to the family \,$\{H_{\lambda}\}$.\, Let, as a family in \,\,$Hilb^{P(n)}_{{\bf P}^3\times{\bf P}^2}(\Delta)$\,\,, it correspond to the map \,$y: \Delta \rightarrow H$. 


Now, the maps $z$\, and $\psi\circ y$\, coincide on $\Delta^*$\,. Since $H$ \,is complete, they coincide on the whole $\Delta$\,. This means that under the morphism $\pi$\, the fibre $\tilde {Y}_0$ \, corresponds to the curve $\tilde {V}_0$\,, and hence is irreducible and isomorphic to the blow up of ${\bf P}^3$\, in that curve. Hence, $\tilde {Y}_0 \cong \tilde {Y}_0'$\,\,, and we conclude that $\tilde {Y}'$ can be realized as a flat family in ${\it Hilb}_{{\bf P}^3\times{\bf P}^2}^{P(n)}(\Delta)$\,\, with a group action $\{g_t\}$\, as above.
\par$ $\par
It is straightforward to verify that the list of all curves in ${\bf P}^3$\, defined by some ideal $I = (\Delta_1, \Delta_2, \Delta_3)$\, contains the following ones (compare with the results in \cite{Piene}):
\begin{enumerate}
\item twisted cubic curves,
\item connected unions of a conic and a line which span the projective space,
\item connected unions of three nonconcurrent lines which span the projective space,
\item the unions of three concurrent lines which span the projective space,
\item unions of a line and a double line given by the ideal $(X_0X_2, X_2^2, X_0X_3)$,
\item the unions of a line and a double line given by the ideal $(X_0^2, X_0X_1, X_0X_2 - X_2X_3)$,
\item triple lines given by the ideal $(X_0X_2 - X_1^2, X_1X_2, X_2^2)$,
\item triple lines given by the ideal $(X_0^2, X_0X_1, X_1^2)$.
\end{enumerate}

Only the last curve corresponds to a reducible intersection of ${\bf P}^3\times{\bf P}^2 $\, with ${\bf P}^9$\, in ${\bf P}^{11}$ - it is an union of the blowup of ${\bf P}^3$\, in the line $ \{X_0 = 0, X_1 = 0\}$\,\, and ${\bf P}^1\times{\bf P}^2$.

Summing up, we have the following result:

\begin{theo}
Consider the degenerations $\tilde{Y} \subset \Delta\times({\bf P}^3\times{\bf P}^2)$\,\, of $V_{38}$\,,\, for which there is a group $\{g_t\}$\, acting on it and lifting the dilation morphisms of $\Delta$\,, and with irreducible central fibre $\tilde{Y}_0$.\, Then, every such degeneration determines a degeneration $\tilde{V} \subset \Delta\times{\bf P}^3$\, of a twisted cubic curve in ${\bf P}^3$\, to one of the curves 1),\ldots,7) above, together with a group $\{h_t\}$\, acting on $\tilde{V}$\, and lifting the dilation morphisms of $\Delta$\,.\, The latter degeneration determines, up to action of $PGL(3)$\,, the pair $(\tilde{Y}, g_t)$\,. If we identify the groups $Aut(\tilde{Y}_0)$\,\, and \,$Aut(\tilde{V}_0 \subset {\bf P}^3)$\,\,, then $h_t = g_t$\, for all $t$.

\end{theo}

\subsection{Toric degenerations of $V_{38}$\, in ${\bf P}^3\times{\bf P}^2$.}

We are going to apply Theorem \ref{torus} and to calculate the Generalized Futaki Invariant of the special toric degenerations of $V_{38}$\, in ${\bf P}^3\times{\bf P}^2$.

Let $X = {\bf P}^9 \cap \, ({\bf P}^3\times{\bf P}^2)$\,\, be an irreducible intersection which is a toric variety, and let denote by $I$\, the ideal of the curve $C \subset {\bf P}^3$\,, such that $X$ \,is the blowup of ${\bf P}^3$\, with respect to  $I$.

Since $Aut(X) \cong Aut(C \subset {\bf P}^3)$\,\,, the only possibilities for $I$\, are:
\begin{enumerate}
\item $I = I_1 = (X_0X_2, X_0X_3, X_2X_3 )$\,, and $C = C_1$\, is a connected union of three nonconcurrent lines in ${\bf P}^3$\, which span ${\bf P}^3$;
\item $I = I_2 = (X_1X_2, X_1X_3, X_2X_3 )$\,, and $C = C_2$\, is an union of three concurrent lines in ${\bf P}^3$\, which span ${\bf P}^3$;
\item $I = I_3 = (X_1X_2, X_1X_3, X_2^2 )$\,, and $C= C_3$\, is an union of the double line with ideal $(X_1, X_2^2 )$\, and the line with ideal $(X_2, X_3 )$.\end{enumerate}

It is immediate that the variety $X_i$\, corresponding to $I_i\,, i=1, 2, 3$\,\, is normal (it locally is a hypersurface, and is smooth in codimension one), and, hence, it is a toric variety.

We will construct a fan of $X_i$\, starting with the following one of ${\bf P}^3$\,: 
$$ \Delta = \{<e_0, e_1, e_2 >, <e_0, e_1, e_3 >, <e_0, e_2, e_3 >, <e_1, e_2, e_3 >\},$$
where $e_0 + e_1 + e_2 + e_3 = 0$\, , and the lattice where $\Delta$\, lives is generated by $\{e_1, e_2, e_3\}$\,. In other words, we consider ${\bf P}^3$\, as a compactification of ${\bf C}^3$\, with the standard (coordinatewise ) action of the torus ${\bf C}^*\times{\bf C}^*\times{\bf C}^*$.

\subsubsection{Description of the case $X = X_1$.}

The variety $X_1$\, has two terminal singular points (above the points $(0: 1: 0: 0)$\, and $(0: 0: 1: 0 )$\, of ${\bf P}^3$)\,, it has Picard number $b_2(X_1) = 7$\,, and a group of holomorphic automorphisms
$$Aut_0(X_1) = \left(\begin{array}{cccc}
* & 0 & 0 & 0\\ {*} & * & 0 & 0\\0 & 0 & * & *\\0 & 0 & 0 & * \end{array}\right)/
{\bf C}^*.$$   

The corresponding fan $\Delta _1$\, consists of the cones
$$<e_1, e_2, e_2 + e_3 >, <e_1, e_3, e_2 + e_3 >, <e_3, e_0 + e_3, e_2 + e_3 >,$$ $$ < e_0, e_2, e_0 + e_1 >,
 <e_0, e_2, e_0 + e_1 >, <e_0, e_0 + e_3, e_0 + e_1 >,$$ $$ < e_0, e_2, e_0 + e_3, e_2 + e_3 >, <e_1, e_3, e_0 + e_1, e_0 + e_3 >.$$

The polytope $-P_{-\omega_X}$\, is the convex hull of the points

$$(1, 1, 0 ), (1, 0, 1 ), (-1, 0, 1 ), (0, 1, -2 ),$$ $$ (1, 1, -2 ), (0, -1, 0 ), (-2, 0, 1 ), (1, -2, 1 ).$$

It has a baricentre
$$ {\bf b}_1 = {1\over 38}.({2\over 3}, {2\over 3}, - {2\over 3}).$$  

\subsubsection{Description of the case $X = X_2$.}

The variety $X_2$\, has three terminal singular points (above the point $(1: 0: 0: 0 )$\, of ${\bf P}^3$\,), has Picard number $b_2(X_2) = 8$\,, and a group of holomorphic automorphisms
$$Aut_0(X_2) = \left(\begin{array}{cccc}
* & * & * & *  \\ 0 & * & 0 & 0  \\ 0 & 0 & * & 0  \\ 0 & 0 & 0 & * \end{array} \right) / {\bf C}^*.$$

The corresponding fan $\Delta _2$\, consists of the cones
$$<e_0, e_1, e_1 + e_2 >, <e_0, e_2, e_1 + e_2 >, <e_0, e_2, e_2 + e_3 >,$$
$$<e_0, e_3, e_2 + e_3 >, <e_0, e_3, e_1 + e_3 >, <e_0, e_1, e_1 + e_3 >,$$
$$<e_1, e_1 + e_2, e_1 + e_3, e_1 + e_2 + e_3 >, <e_2, e_1 + e_2, e_2 + e_3, e_1 + e_2 + e_3 >.$$

The polytope $- P_{-\omega_{X_2}}$\, is the convex hull of the points
$$(1, 0, 0 ), (0, 1, 0 ), (0, 0, 1 ), (1, 0, -2 ), (0, 1, -2 ),$$
$$(1, -2, 0 ), (-2, 1, 0 ), (-2, 0, 1 ), (0, -2, 1 ).$$

It has a barycentre
$${\bf b}_2 = {1\over 38}.(-{2\over 3}, -{2\over 3}, -{2\over 3}).$$ 

\subsubsection{Description of the case $X = X_3$.}

The variety $X_3$\, has three $T$-invariant singular points: $Q_1, Q_2, Q_3$.\, One of them, say $Q_1$\, is a terminal singular point, while $Q_2$ is a quotient singularity, and $Q_3$ is a canonical singularity. $Q_2$\, and $Q_3$ lie on a $T$-invariant ${\bf P}^1$\, which consists of terminal singular points (outside $Q_2$ and $Q_3$). The Picard number $b_2(X_3) = 8$,\, and the group of holomorphic automorphisms is
$$Aut_0(X_3) = \left(\begin{array}{cccc}
* & * & * & *  \\ 0 & * & 0 & 0  \\0 & 0 & * & 0  \\0 & 0 & * & *\end{array} \right)/{\bf C}^*.$$
The fan $\Delta_3$\, consists of the cones
$$<e_0, e_1, e_3 >, <e_0, e_2, e_2 + e_3 >, <e_0, e_3, e_2 + e_3 >, $$
$$<e_0, e_1, 2e_1 + e_2 >, <e_0, e_2, 2e_1 + e_2 >, <-e_0, e_3, e_2 + e_3 >,$$
$$<-e_0, e_1, e_3, 2e_1 + e_2 >, <-e_0, e_2, e_2 + e_3, 2e_1 + e_2 >.$$
The polytope $-P_{-\omega_{X_3}}$\, is the convex hull of the points
$$(1, -3, 1 ), (-2, 1, 0 ), -2, 0, 1 ), (1, -1, -1),$$
$$(0, 1, -2 ), (0, 0, 1 ), (1, -1, 1 ), (0, 1, 0 ).$$
It has a barycentre
$${\bf b}_3 = {1\over 38}.(-{2\over 3}, -2, {2\over 3}).$$    

By a direct computation one can show that for each special degeneration of \,$V_{38}$\, to one of the toric varieties $X_i$\,, the Generalized Futaki Invariant of \,$X_i$\, evaluated on the admissible vector field, corresponding to that degeneration, has positive real part ($i = 1, 2, 3.$).

\end{document}